\documentclass[11pt,a4paper,reqno]{amsart}
\usepackage{amsfonts,amsthm,latexsym,amsmath,amssymb,amscd,amsmath,epsf,morefloats}
\usepackage{placeins}
\usepackage{graphicx}
\usepackage{tikz}
\usetikzlibrary{arrows}
\usetikzlibrary{patterns}

\newtheorem{theo+}           {Theorem}
\newtheorem{prop+}           {Proposition}
\newtheorem{coro+}           {Corollary}
\newtheorem{lemm+}           {Lemma}

\theoremstyle{definition}
\newtheorem{defi+}           {Definition}
\newtheorem{problem}         {Problem}

\theoremstyle{remark}
\newtheorem{rema+}           {Remark}

\newenvironment{theorem}{\begin{theo+}}{\end{theo+}}

\newenvironment{lemma}{\begin{lemm+}}{\end{lemm+}}

\usepackage{float}
\usepackage[caption = false]{subfig}

\title[Temperature on rods with Robin boundary conditions  ]{Temperature on rods with Robin boundary conditions}
\author{Dimitrios Betsakos and Alexander~Solynin}
\date{\today}

\date{\today}


\begin{document}
         \numberwithin{equation}{section}

\begin{abstract} %
 We consider solutions $u_f$ to the one-dimensional Robin problem with the heat source $f\in L^1[-\pi,\pi]$
 and Robin parameter $\alpha>0$. For given $m$, $M$, and $s$,
 $0\le m<s<M$, we identify the heat sources
 $f_0$, such that $u_{f_0}$ maximizes the temperature gap
 $\max_{[-\pi,\pi]}u_f -\min_{[-\pi,\pi]}u_f$ over all heat
 sources $f$ such that $m\le f\le M$ and $\|f\|_{L^1}=2\pi s$. In
 particular, this answers a question raised by J.~J.~Langford and
P.~McDonald in \cite{LM}. We also identify heat sources, which
maximize/minimize $u_f$ at a given point $x_0\in [-\pi,\pi]$ over
the same class of heat sources as above and discuss a few related questions. %
\end{abstract} %

\maketitle


 \section{Heating with Robin boundary
conditions.} %
Recently,  J.~J.~Langford and
P.~McDonald \cite{LM} studied the one-dimensional Poisson equation with Robin boundary conditions. They considered the following physical setup:
Suppose that a metal rod of length $2\pi$ is located along the interval $[-\pi,\pi]$. Suppose that to half of the locations of the rod, heat is generated uniformly; call this set $E$.  On the remaining half, heat is neither generated nor absorbed. The ends of the rod interact with the cooler environment so that there is a heat flux from the rod which is proportional to the temperature at each end (Newton's law of cooling).

Let $u$ be the steady-state temperature function; it satisfies the Poisson equation
\begin{equation}\label{0.1}
-u''(x)=\chi_E(x),\;\;\;x\in [-\pi,\pi]
\end{equation}
with Robin boundary conditions
\begin{equation}\label{0.2}
-u'(-\pi)+\alpha u(-\pi)=u'(\pi)+\alpha u(\pi)=0,
\end{equation}
where $\alpha>0$ and  $\chi_E$
stands for the characteristic function of the set $E$.

Langford and
McDonald \cite{LM} studied the problem of where one should locate the heat sources to maximize the hottest steady-state temperature. In other words, for which set $E$, for the solution $u$ of the boundary value problem (\ref{0.1})-(\ref{0.2}), $\max_{[-\pi,\pi]} u$ is maximal.
They showed that this quantity is maximal when $E$ is an {\it interval located symmetrically} in the middle of the rod; namely $E=[-\pi/2,\pi/2]$.
The authors of \cite{LM} actually studied a more general problem and obtained much stronger comparison results. They observed, however, that the symmetric interval is not extremal for another problem.  They considered the temperature gap
over $[-\pi,\pi]$, i.e. the quantity %
$$ 
{\rm osc}(u):=\max_{[-\pi,\pi]} u-\min_{[-\pi,\pi]} u %
$$  
and showed that ${\rm osc}(u)$ is not maximized for
$E=[-\pi/2,\pi/2]$. So they raised the following question:

\begin{problem}
 Where should we place the heat
sources to maximize the temperature gap?
\end{problem}

They suggested that the extremal set $E$ is again an interval which, however, is not symmetrically located on $[-\pi,\pi]$. In the present note, we will study this conjecture. As in \cite{LM}, we will consider a more general setting.

\medskip

The Robin problem for $f\in L^1[-\pi,\pi]$ and $\alpha>0$ is to
find
$u\in C^1[-\pi,\pi]$ such that %
\begin{enumerate}  %
\item[1.] $u'$ is absolutely continuous on $[-\pi,\pi]$, %
\item[2.] $-u''=f$ a.e. on $(-\pi,\pi)$, %
\item[3.] $-u'(-\pi)+\alpha u(-\pi)=u'(\pi)+\alpha u(\pi)=0$. %
\end{enumerate} %

It was shown in Proposition~2.1 in \cite{LM} that the Robin
problem has a unique solution given by the equation %
\begin{equation}\label{1.1}
u_f(x)=\int_{-\pi}^\pi G(x,y)f(y)\,dy.
\end{equation}
Here, $G(x,y)$ stands for the Green's function for Robin problem,
which is %
\begin{equation}\label{1.2}
G(x,y)=-\frac{1}{2}c_\alpha
xy-\frac{1}{2}|x-y|+\frac{1}{2c_\alpha},\quad x,y\in [-\pi,\pi],
\end{equation}
where %
\begin{equation}\label{1.3} %
c_\alpha=\frac{\alpha}{1+\alpha\pi}. %
\end{equation} %

To recall a few basic facts about solutions of the Robin problem, we
note first that, as simple Calculus shows,  $G(x,y)>0$ for all
$x,y\in [-\pi,\pi]$. Another simple but important conclusion from
(\ref{1.1}) is that the solution to the Robin problem is an
additive function of the heat source; i.e. if $m_1$, $m_2$ are
constants and $f_1,f_2\in
L^1$ (here and below $L^1$ stands for $L^1[-\pi,\pi]$), then%
\begin{equation}\label{1.4} %
 u_f=m_1u_{f_1}+m_2u_{f_2}, \quad {\mbox{where $f=m_1f_1+m_2f_2$.}}%
\end{equation} %

Furthermore, if $f\ge 0$, then it is immediate from property~2
above that $u_f$ is concave on $[-\pi,\pi]$ and therefore it takes
its minimal value at one of the end points $x=\pm\pi$ and it takes
its maximal value either at a single point or on some closed
subinterval of $[-\pi,\pi]$. If $u_f$ takes its minimal value at
$-\pi$, it follows from property~3 that $u_f(x)\ge
u_f(-\pi)=\alpha^{-1}u'_f(-\pi)>0$ and therefore, in this case,
$u_f(x)>0$ for $x\in [-\pi,\pi]$ unless $f\equiv 0$. The same
conclusion follows if $u_f$ takes its minimal value at $\pi$.

  To put the
question in Problem~1 in a more general setting,  we consider, for
given $m$, $M$, and $s$ such that $0\le m<s< M$, the class
$\mathcal{F}=\mathcal{F}(m,M,s)$ of heat sources $f\in
L^1$ such that  $m\le f(x)\le M$ for $x\in [-\pi,\pi]$
and $\|f\|_{L^1}=2\pi s$. The parameters $m$, $M$, and $s$ can be
interpreted as the ground heat, the top heat, and the average heat
over the rod.

Our main goal in this note is to prove the following theorem,
which solves the maximal temperature gap problem for the class
$\mathcal{F}(m,M,s)$ and therefore, as a special case, it provides
a solution to Problem~1. In this theorem and below,
$I(a,l)\subset \mathbb{R}$ denotes the closed interval of length
$2l>0$ centered at $a$.
\begin{theorem}\label{Theorem 1}%
Let $u_f$ solves the Robin problem for $f\in \mathcal{F}(m,M,s)$ and $\alpha>0$. Then %
\begin{equation}\label{1.5}
{\rm osc}(u_f)\leq (M-m)\Theta_\alpha(l,\delta),
\end{equation}
where $l=\pi(s- m)/(M-m)$, $\delta=m/(M-m)$  and
$\Theta_\alpha(l,\delta)$ is defined by equations (\ref{4.10}),
(\ref{4.11}) and (\ref{4.12}) in Section~4.

Equality holds in (\ref{1.5}) if and only if $f(x)=f_0(x)$ or
$f(x)=f_0(-x)$  a.e. on $[-\pi,\pi]$, where
$f_0=m+(M-m)\,\chi_{I(a_g,l)}$ with $l=\pi(s- m)/(M-m)$  and $a_g$
defined by equation (\ref{4.9}) in Section~4.
\end{theorem} %

For $m=1$, $M=3$, and $s=\frac{7}{5}$, the graph of the maximal
temperature gap  $(M-m)\Theta_\alpha(l,\delta)$ considered as a
function of $\alpha$ is shown in Figure~1(a). Figure~1(b) displays
the graph of the extremal function $u_{f_0}(x)$, where
$f_0=m+(M-m)\,\chi_{I(a_g,l)}$ with $m=1$, $M=3$, $s=\frac{7}{5}$
and $\alpha=1/2$.

\begin{figure}
\subfloat[ ]{\includegraphics[width = 210pt]{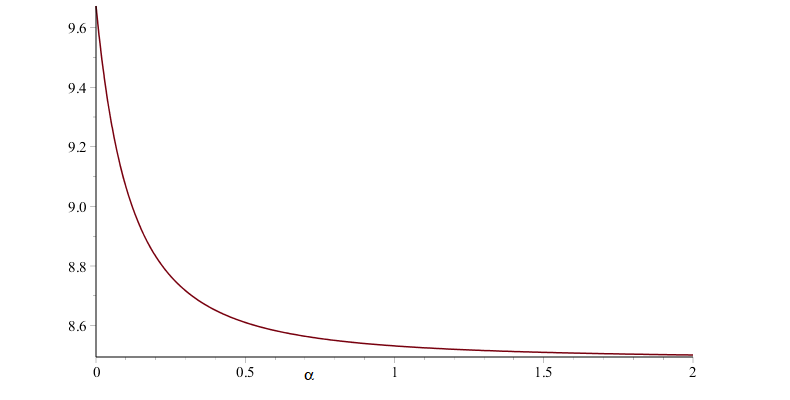}}
\subfloat[ ]{\includegraphics[width = 145pt]{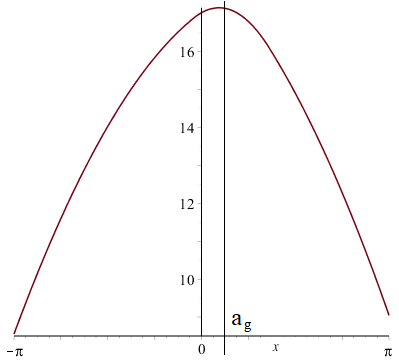}}
\caption{(a) The maximal temperature gap
$(M-m)\Theta_\alpha(l,\delta)$ for $m=1$, $M=3$, $s=\frac{7}{5}$
and $0<\alpha<2$; (b) Extremal function $u_{f_0}(x)$ for $m=1$,
$M=3$, $s=\frac{7}{5}$ and $\alpha=1/2$.} \label{some example}
\end{figure}

Returning to the context of Problem 1, let us suppose that $E$ is
a measurable subset of $[-\pi,\pi]$ of length (one-dimensinal
Lebesgue measure) equal to $\pi$. We apply Theorem \ref{Theorem 1}
with $f=\chi_E$, $m=0$, $M=1$, $\delta=0$, and $l=\pi/2$. The
parameter $\alpha_0$ defined by equation (\ref{4.7.1}) in Section
4 takes the value $\alpha_0=\frac{2}{\sqrt{3}\pi}$. Let
$$
E^*=I(a_g,\pi/2)=[a_g-\pi/2,a_g+\pi/2]\;\;\;\;\hbox{and}\;\;\;\;f^*=\chi_{E^*},
$$
where (see formula (\ref{4.9}))
$$
a_g=\begin{cases}\pi/2,
\;\;\;\;\;\;\;\;\;\;\;\;\;\;\;\;\;\;\;\;\;\;\;\;\;\;\;\hbox{if}\;\;0<\alpha\leq
\frac{2}{\sqrt{3}\pi},\\
\frac{\pi/2}{(1+\alpha\pi)(\pi c_\alpha/2-\pi^2 c_\alpha^2/4)},
\;\;\hbox{if}\;\;\alpha\ge \frac{2}{\sqrt{3}\pi}.
\end{cases}
$$
It follows (cf. \cite[Proposition 3.3]{LM}) that when
$0<\alpha\leq \frac{2}{\sqrt{3}\pi}$, we have $E^*=[0,\pi]$; when
$\alpha\ge \frac{2}{\sqrt{3}\pi}$, the location of the interval
$E^*$ depends on $\alpha$ and as $\alpha$ increases, $E^*$ moves
from the right end to the center. By Theorem \ref{Theorem 1}, we
have ${\rm osc}(u_f)\leq  {\rm osc}(u_{f^*})$. The solution to
Problem 1 is: The temperature gap is maximized uniquely when the
heat sources are placed on $E^*$ or on $-E^*$.


\medskip

Another interesting problem on the distribution of heat on a rod
is to identify heat sources $f\in \mathcal{F}(m,M,s)$, which
generate the maximal possible temperature and the minimal possible
temperature at a fixed location $x_0\in [-\pi,\pi]$ of the rod.
Notice that if $f^-(x)=f(-x)$, then $u_{f^-}(x)=u_f(-x)$. Thus,
working with this problem, we may assume that $x_0\in [0,\pi]$.
Its solution is given by the following theorem.

\begin{theorem}  \label{Theorem 2}%
Let $u_f$ solves the Robin problem for $f\in \mathcal{F}(m,M,s)$ and $\alpha>0$ and let $x_0\in [0,\pi]$ be fixed. Then %
\begin{equation}\label{1.6}
\eta_\alpha(x_0)M-(M-m)\nu_\alpha(x_0,l^-)\le u_f(x_0)\le
\eta_\alpha(x_0)m+(M-m)\nu_\alpha(x_0,l),
\end{equation}
where $l=\pi(s- m)/(M-m)$, $l^-=\pi-l$, %
and the functions $\eta_\alpha(x)$ and $\nu_\alpha(x,l)$ are
defined in Section~2 by equations (\ref{2.2}) and (\ref{2.4}),
respectively.

Equality holds in the right inequality in (\ref{1.6}) if and only
if $f=f_0^+$ a.e. on $[-\pi,\pi]$, where
$f_0^+=m+(M-m)\,\chi_{I(a_m,l)}$
 with $l$ defined above and $a_m=x_0(1-lc_\alpha)$ if $x_0(1-lc\alpha)<\pi-l$ and $a_m=\pi-l$ otherwise. %

Equality holds in the left inequality in (\ref{1.6}) if and only
if $f=f_0^-$ a.e. on $[-\pi,\pi]$, where
$f_0^-=M-(M-m)\,\chi_{I(a_m^-,l^-)}$ with $l^-=\pi-l$
and $a_m^-=x_0(1-l^-c_\alpha)$ if $x_0(1-l^-c\alpha)<\pi-l^-$ and $a_m^-=\pi-l^-$ otherwise. %
\end{theorem} %





It would be also useful to know how warmer a fixed spot $x_0$
could be compared to the edges of the rod. The answer to this
question is the following.

\begin{theorem}  \label{Theorem 3}%
Let $u_f$ solves the Robin problem for $f\in \mathcal{F}(m,M,s)$ and $\alpha>0$ and let $x_0\in [-\pi,\pi]$ be fixed. Then %
\begin{equation}\label{1.7}
u_f(x_0)-u_f(-\pi)\le m
(\eta_\alpha(x_0)-\eta_\alpha(-\pi))+(M-m)\tau_\alpha(x_0,l),
\end{equation}
where $l=\pi(s- m)/(M-m)$, %
the function $\eta_\alpha(x)$ is defined by equation (\ref{2.2})
and the function $\tau_\alpha(x,l)$ is defined by equation
(\ref{3.5}).

Equality holds in (\ref{1.7}) if and only if $f=f_e$ a.e. on
$[-\pi,\pi]$, where $f_e=m+(M-m)\,\chi_{I(a_e,l)}$
 with $l$ defined above and $a_e=a_e(x,l)$ defined in equation (\ref{3.4}). %
\end{theorem} %

The main results of \cite{LM} stated in Theorems 1.3 and 1.5
concern the comparison principles for heating problems with the
Robin  and Newmann boundary conditions, respectively. In these
problems, the extremal distribution of heat is symmetric with
respect to the center of the rod. With this symmetry, the authors
of \cite{LM} were able to use symmetrization methods due to
G.~Talenti \cite{T} and A.~Baernstein~II \cite{B} to prove their
theorems. We want to mention here that S.~Abramovich in her paper
\cite{A} published in 1975 already used the symmetrization method
due to G.~P\'{o}lya and G.~Szeg\"{o} \cite{PS} to prove
interesting results on the eigenvalues of the differential system
$y''(x)+\lambda p(x) y(x)=0$, $y(\pm 1)=0$ for the function
$p(x)\ge 0$ defined on the string $(-1,1)$. More recently, the
symmetrization method similar to the one used in \cite{LM} in
combination with the polarization technique was used in \cite{BS}
to study several problems on heat distribution in the cylindrical
pipes heated along various regions on the surface area.

We stress here that the extremal distributions of heat in our
Theorems~1, 2 and 3 are not symmetric with respect to the center
of the rod, in general. Thus, the classical symmetrization
technique cannot be applied in these problems while certain
versions of polarization technique used in \cite{BS} still can be
applied.

 \section{Heating a fixed spot by a single interval.} %
 Suppose that 
 the Robin rod $[-\pi,\pi]$ is heated with unit density along the interval $I=I(a,l)$
 centered at the point $a\in (-\pi,\pi)$ with length $2l$ such
 that $-\pi\le a-l<a+l\le \pi$. Thus, we assume here that $f=\chi_I$. Using the integral
 representation (\ref{1.2}) for the Robin temperature with the Green's function given by (\ref{1.3}),
 we evaluate $u_{\chi_I}=u_{\chi_{I(a,l)}}(x,\alpha)$ as follows: %

\begin{eqnarray}\label{2.1}%
u_{\chi_I}&=&\int_I [(-1/2)c_\alpha)xy-(1/2)|x-y|+1/(2c_\alpha)]\,dy \nonumber\\
&{=}& \left\{\begin{array}{lr}
        l[(1-ac_\alpha)x+(c_\alpha^{-1}-a)], &  -\pi\le x\le a-l,\\  
        -\frac{1}{2}x^2+a(1-lc_\alpha)x+\frac{l}{c_\alpha}-\frac{a^2+l^2}{2}, &  a-l<x<a+l,\\   
        l[-(1+ac_\alpha)x+(c_\alpha^{-1}+a)], &  a+l\le x\le \pi.   
        \end{array}\right.  
\end{eqnarray}

In particular, if the whole rod $[-\pi,\pi]$ is heated with unit
density, then $a=0$, $l=\pi$, and
$u_{\chi_{[-\pi,\pi]}}(x)=\eta_\alpha(x)$, where 
\begin{equation}\label{2.2}
\eta_\alpha(x)=
-\frac{1}{2}x^2+\frac{\pi}{\alpha}+\frac{\pi^2}{2},
\quad -\pi\le x\le \pi.
\end{equation}

Next, we fix $x_0\in [0,\pi]$, $l\in (0,\pi)$, $\alpha>0$ and
treat $u_{\chi_I}$ as a function $F(a)$ of the variable $a\in
[-\pi+l,\pi-l]$.

We have to consider the following cases: %
\begin{enumerate}
\item[1)] If $0\le x_0\le -\pi+2l$, then 
\begin{equation}
F(a)=-\frac{1}{2}a^2+x_0(1-lc_\alpha)a-\frac{1}{2}x_0^2+\frac{l}{c_\alpha}-\frac{l^2}{2},
\quad  -\pi+l\le a\le \pi-l, \nonumber
\end{equation}
\item[2)] If $-\pi+2l< x_0< \pi-2l$, then %
\begin{equation}
F(a)=\left\{\begin{array}{ll}
 l(1-x_0c_\alpha)a+l(\frac{1}{c_\alpha}-x_0), &  -\pi+l \le a\le x_0-l,\\ %
 -\frac{1}{2}a^2+x_0(1-lc_\alpha)a-\frac{1}{2}x_0^2+\frac{l}{c_\alpha}-\frac{l^2}{2},
 & x_0-l\le a\le x_0+l,\\
  -l(1+x_0c_\alpha)a+l(\frac{1}{c_\alpha}+x_0), &  x_0+l\le a\le \pi-l.
\end{array}\right. \nonumber %
\end{equation} %
\item[3)] If $\max\{-\pi+2l,\pi-2l\}\le x_0\le \pi$, then %
\begin{equation}
F(a)=\left\{\begin{array}{ll}
 l(1-x_0c_\alpha)a+l(\frac{1}{c_\alpha}-x_0), &  -\pi+l \le a\le x_0-l,\\ %
 -\frac{1}{2}a^2+x_0(1-lc_\alpha)a-\frac{1}{2}x_0^2+\frac{l}{c_\alpha}-\frac{l^2}{2},
 & x_0-l\le a\le \pi-l.
\end{array}\right. \nonumber %
\end{equation} %
\end{enumerate} %

A simple argument, left to the interested reader, shows that, in
all three cases, $F(a)$ is positive and concave on the interval
$[-\pi+l,\pi-l]$ and takes its minimal value
$\mu_\alpha=\mu_\alpha(x_0,l)$ at
$a=-\pi+l$. Evaluating $\mu_\alpha=F(-\pi+l)$, we find %
\begin{equation}
\mu_\alpha=\left\{\begin{array}{l}
 l[((\pi-l)c_\alpha-1)x_0+\frac{1}{c_\alpha}+l-\pi],  {\mbox{\  if $-\pi+2l<x_0\le \pi$,}} \\
  -\frac{1}{2}x_0^2-(\pi-l)(1-lc_\alpha)x_0-l^2+l(\frac{1}{c_\alpha}+\pi)-\frac{\pi^2}{2},  {\ \mbox{otherwise.}} %
\end{array}\right. \nonumber %
\end{equation} %

Next, we find the maximum $\nu_\alpha=\nu_\alpha(x_0,l)=\max F(a)$
taken over the interval $-\pi+l\le a\le \pi-l$ and identify the
point $a_m\in [-\pi+l,\pi-l]$, where this maximum is achieved. Let
$q(a)$ denote the quadratic function as in parts 1)--3). Then
$q(a)$ takes its maximum at the point $a_0=x_0(1-lc_\alpha)$.
Notice that $0\le a_0\le x_0$ with equality sign in either of
these inequalities if and only if $x_0=0$. If $a_0+l\le \pi$, then
the function $F(a)$ is defined at $a_0$ and takes its maximum at
this point. Thus, if $a_0+l\le \pi$, then $a_m=x_0(1-lc_\alpha)$.
If $a_0+l> \pi$, then $F(a)$ takes its maximum at $a_m=\pi-l$.
Combining these cases, we have the following equation for the
central point $a_m=a_m(x_0,l,\alpha)$ of the heating interval of
length $2l$, which
generates the maximal temperature at the point $x_0$: %
\begin{equation}\label{2.3}
a_m=\left\{\begin{array}{ll}
x_0(1-lc_\alpha), & {\mbox{if  $x_0<\frac{\pi-l}{1-lc_\alpha}$,}}\\ %
 \pi-l, & {\mbox{otherwise.}}
\end{array}\right.
\end{equation} %

 With these
notations, we can evaluate the maximum
$\nu_\alpha=\nu_\alpha(x_0,l)$ as follows:

\begin{equation}\label{2.4}
\nu_\alpha=\left\{\begin{array}{lr}
\frac{l}{c_\alpha}\left(1-\frac{l
c_\alpha}{2}\right)\,\left(1-x_0^2 c_\alpha^2\right), & {\mbox{if
$x_0<\frac{\pi-l}{1-lc_\alpha}$,}}
\\ %
\frac{1}{2}(\pi-x_0)^2+l(1-c_\alpha x_0)(2\pi-l+\frac{1}{\alpha}),
&  {\mbox{otherwise.}}
\end{array} \right.
\end{equation}

The inequality $a_0+l\le \pi$ is
equivalent to the inequality %
$$ 
\alpha\ge \frac{x_0+l-\pi}{(\pi-l)(\pi-x_0)}.
$$ 
Let us define $\alpha_m\ge 0$ as %
\begin{equation}\label{2.5}
\alpha_m=\left\{\begin{array}{ll}
\frac{x_0+l-\pi}{(\pi-l)(\pi-x_0)}, & {\mbox{if  $x_0<\pi-l$,}}\\ %
 0, & {\mbox{otherwise.}}
\end{array}\right.
\end{equation} %
Now, our arguments above show that if $x_0\in [0,\pi)$, $l\in
(0,\pi)$ and $\alpha>\alpha_m$, then  the function $F(a)$ achieves
its maximum $\nu_\alpha(x_0,l)$,  given by the first line of
(\ref{2.4}), at the point $a_m=x_0(1-lc_\alpha)$, and if
$0<\alpha\le \alpha_m$, then $F(a)$ achieves its maximum
$\nu_\alpha(x_0,l)$, given by the second line of  (\ref{2.4}), at
the point $a_m=\pi-l$.

Combining our results and using the notation introduced above, we
obtain the following lemma.

\begin{lemma}  %
{\rm 1)} Let $x_0\in [0,\pi]$, $l\in (0,\pi)$, and $\alpha>0$ be
fixed and let $a$ varies from $-\pi+l$ to $\pi-l$.

If $\alpha>\alpha_m$, then the function
$F(a)=u_{\chi_{I(a,l)}}(x_0,\alpha)$ increases from its minimal
value $\mu_\alpha(x_0,l)$ to its maximal value $\nu_\alpha(x_0,l)$
as $a$ varies from $-\pi+l$ to $a_m=x_0(1-lc_\alpha)<\pi-l$, and
$F(a)$ decreases as $a$ varies from $a_m$ to $\pi-l$.

If $0<\alpha\le\alpha_m$, then the function
$F(a)=u_{\chi_{I(a,l)}}(x_0,\alpha)$ increases from
$\mu_\alpha(x_0,l)$ to $\nu_\alpha(x_0,l)$ as $a$ varies from
$-\pi+l$ to $a_m=\pi-l$.

{\rm 2)} Furthermore, the point $a_m$, where $F(a)$ takes its
maximum, stays at $\pi-l$ for $0<\alpha\le \alpha_m$ and $a_m$
decreases from $\pi-l$ to $\frac{\pi-l}{\pi}x_0$, when $\alpha$
runs from $\alpha_m$ to $\infty$.

{\rm 3)} Moreover, if $x_0\in [0,\pi]$ and $\alpha>0$ are fixed
and $a_m$ is considered as a function $a_m(l)$ of $l$, then if
$0<l_1<l_2<\pi$, then
$$
x_0\in
[a_m(l_1)-l_1,a_m(l_1)+l_1]\subset [a_m(l_2)-l_2,a_m(l_2)+l_2]. %
$$
\end{lemma} %


 \section{Temperature gap between a fixed spot and the edges of a rod for a single
 interval.}  %
 As in the previous section, we assume that the Robin rod $[-\pi,\pi]$ is heated
 with unit density along the interval $I=I(a,l)$.
 Let us fix  $x_0\in
 [-\pi,\pi]$, $l\in(0,\pi)$, $\alpha>0$, and consider the temperature gap between
 the point $x_0$ and the left edge of the rod as a function of $a\in
 [-\pi+l,\pi-l]$; i.e. we consider the function %
 \begin{equation} \label{3.1} %
E(a)=u_{\chi_I}(x_0)-u_{\chi_I}(-\pi).
 \end{equation} %
 To find $E(a)$, we use equation (\ref{2.1}). Depending on the
values of $l$ and $x_0$, we have to consider the following cases:
\begin{enumerate}
\item[1)] If $-\pi\le x_0\le \min\{-\pi+2l,\pi-2l\}$, then %
\begin{equation}
E(a)=\left\{\begin{array}{ll}
-\frac{1}{2}a^2+(x_0(1-lc_\alpha)+l(1-\pi c_\alpha))a-
\\ -\frac{1}{2}x_0^2+\pi l-\frac{l^2}{2},
 & {\mbox{if $-\pi+l\le a\le x_0+l$,}}\\
  -lc_\alpha(\pi+x_0)a+l(\pi+x_0), &  {\mbox{if $x_0+l\le a\le
  \pi-l$.}}
\end{array}\right. \nonumber %
\end{equation} %
\item[2)] If $-\pi+2l< x_0< \pi-2l$, then %
\begin{equation}
E(a)=\left\{\begin{array}{ll}
 l(2-(\pi+x_0)c_\alpha)a+l(\pi-x_0), &  {\mbox{if $-\pi+l \le a\le x_0-l$,}}\\ %
-\frac{1}{2}a^2+(x_0(1-lc_\alpha)+l(1-\pi c_\alpha))a-
\\
-\frac{1}{2}x_0^2+\pi l-\frac{l^2}{2},
 & {\mbox{if $x_0-l\le a\le x_0+l$,}}\\
  -lc_\alpha(\pi+x_0)a+l(\pi+x_0), &  {\mbox{if $x_0+l\le a\le
  \pi-l$.}}
\end{array}\right. \nonumber %
\end{equation} %
\item[3)] If $\max\{-\pi+2l,\pi-2l\}\le x_0\le \pi$, then %
\begin{equation}
E(a)=\left\{\begin{array}{ll}
  l(2-(\pi+x_0)c_\alpha)a+l(\pi-x_0), & {\mbox{if $-\pi+l \le a\le x_0-l$,}}\\  %
 -\frac{1}{2}a^2+(x_0(1-lc_\alpha)+l(1-\pi
c_\alpha))a- \\ %
-\frac{1}{2}x_0^2+\pi l-\frac{l^2}{2},
 & {\mbox{if $x_0-l\le a\le \pi-l$.}}
\end{array}\right. \nonumber %
\end{equation} %
\item[4)] If $\pi-2l \le x_0\le -\pi+2l$, then for $-\pi+l\le
a\le \pi-l$, 
\begin{equation}
E(a)=-\frac{1}{2}a^2+(x_0(1-lc_\alpha)+l(1-\pi
c_\alpha))a-\frac{1}{2}x_0^2+\pi l-\frac{l^2}{2}. \nonumber
\end{equation}
\end{enumerate} %

Next, we show that in each of the four cases above there is a
unique point $a_e\in [-\pi+l,\pi-l]$, where the function $E(a)$
achieves its maximum, we call it $\tau_\alpha=\tau_\alpha(x_0,l)$.

On can easily check that in all cases, $E'(-\pi+l)>0$ and, in the
cases 1) and 2), $E'(\pi-l)<0$. Since $E(a)$ is a smooth at most
quadratic function, we conclude from this that, in the cases 1)
and 2), $E(a)$ achieves its maximum $\tau_\alpha$ at the point
$a_e=x_0+l-lc_\alpha(\pi+x_0)$, $-\pi+l<a_e<x_0+l$.

In the cases 3) and 4), we find that %
\begin{equation} \label{3.2} %
E'(\pi-l)=-\pi+2l+x_0-l(\pi+x_0)c_\alpha.
\end{equation} %
Since, in the cases 3) and 4),  $-\pi+2l+x_0>0$, we conclude from
(\ref{3.4}) and (\ref{1.3}) that $E'(\pi-l)<0$ if
$\alpha>\alpha_e$ and $E'(\pi-l)\ge 0$ if $0< \alpha\le \alpha_e$,
where $\alpha_e$ is defined as%
\begin{equation}  \label{3.3} %
\alpha_e=\left\{\begin{array}{ll}\frac{-\pi+2l+x_0}{(\pi-l)(\pi-x_0)}
& {\mbox{if $x_0>\pi-2l$,}} \\
0 & {\mbox{if $x_0\le \pi-2l$.}} %
\end{array} \right.
\end{equation} %
In the first of these cases, $E(a)$ takes its maximum
$\tau_\alpha$
at the point $a_e=x_0+l-lc_\alpha(\pi+x_0)$. 
In the second case, 
 $E(a)$ achieves its maximum  at the point $a_e=\pi-l$.
%

Combining these cases, we conclude that $E(a)$ achieves its maximum at the point
\begin{equation}\label{3.4}
a_e=\left\{\begin{array}{lr} %
\pi-l, & {\mbox{if $x_0>\pi-2l$ and $0<\alpha\le \alpha_e$,}} \\
x_0+l-lc_\alpha(\pi+x_0), & {\mbox{otherwise.}}
\end{array}\right.
\end{equation} %

Calculating the values of $\tau_\alpha$ in the cases mentioned
above, we obtain: %
\begin{equation} \label{3.5} %
\tau_\alpha=\left\{\begin{array}{l} %
-\frac{1}{2}x_0^2+\frac{\pi^2}{2}+(x_0(1-lc_\alpha)+l(1-\pi
c_\alpha))(\pi-l), \\ %
{\mbox{\ \ \ \ \ \ \ \ \ \ \ \ \ \ \ \ \ \ \ \ \ \ \ \ \ \ \ \ \ \ \ \ \ \  if $x_0>\pi-2l$ and $0<\alpha\le \alpha_0$,}}\\
\frac{1}{2}(x_0+l-lc_\alpha(\pi+x_0))^2-\frac{1}{2}x_0^2+\pi
l-\frac{1}{2}l^2, {\ \ \ \ \ \mbox{otherwise.}}

\end{array} \right.
\end{equation} %

Combining our results, we obtain the following lemma.

\begin{lemma}  %
Let $x_0\in [-\pi,\pi]$, $l\in (0,\pi)$, and $\alpha>0$ be fixed
and let $a$ vary from $-\pi+l$ to $\pi-l$.

If $x_0$ and $l$ satisfy the inequalities in parts 1) or 2) given
above in this section or if $x_0$ and $l$ satisfy inequalities
given in parts 3) or 4) and also $\alpha>\alpha_e$, where
$\alpha_e$ is defined in (\ref{3.3}), then the function
$E(a)=u_{\chi_I}(x_0)-u_{\chi_I}(-\pi)$ increases from $E(-\pi+l)$
to its maximal value $\tau_\alpha=\tau_\alpha(x_0,l)$ given by the
second line in the equation (\ref{3.5}) as $a$ varies from
$-\pi+l$ to $a_e$ defined by equation (\ref{3.4}), and $E(a)$
decreases as $a$ varies from $a_e$ to $\pi-l$.

If $x_0$ and $l$ satisfy the inequalities given in parts 3) or 4)
and also  $0<\alpha\le\alpha_e$, then the function $E(a)$
increases from $E(-\pi+l)$ to its maximal value $\tau_\alpha$
given by the first line in the equation (\ref{3.5}) as $a$ varies
from $-\pi+l$ to $\pi-l$.
\end{lemma} %


 \section{Maximal temperature gap for a single interval.} %
The goal here is to evaluate the temperature gap ${\rm{osc}}(u_J)$
for the heat source $J=\delta+\chi_{I(a,l)}$, $l\in[0,\pi)$, $a\in
[0,\pi-l]$,
$\delta\ge 0$. It follows from our calculations in Section~2 that %
\begin{equation}  \label{4.1} %
u_J(x)=\delta\eta_\alpha(x)+u_{\chi_{I(a,l)}}(x),
\end{equation} %
where $\eta_\alpha(x)$ is defined by (\ref{2.2}) and
$u_{\chi_{I(a,l)}}$ is defined by (\ref{2.1}).

An elementary calculation shows that, under our assumptions,
$u_J(-\pi)\le u_J(\pi)$ with equality sign if and only if $a=0$.
Since $u_J$ is concave on $[-\pi,\pi]$ it follows that
\begin{equation} \label{4.2}  %
\min_{[-\pi,\pi]} u_J=u_J(-\pi)=\alpha^{-1}(\pi
\delta+l(1-ac_\alpha)).
\end{equation} %

To find the maximum of $u_J$, we note that $u_J$ is a concave
piece-wise at most quadratic function, which achieves its maximum
at the
point %
\begin{equation}  \label{4.3} %
x_0=\frac{a(1-lc_\alpha)}{1+\delta}.
\end{equation}  %
Thus, $0\le x_0\le a$ with equality sign if and only if $a=0$.
Evaluating $u_J(x_0)$ and simplifying, we find: %
\begin{equation}\label{4.4}
\max_{[-\pi,\pi]} u_J=u_J(x_0)
=\frac{1}{2}\left(\frac{(1-lc_\alpha)^2}{1+\delta}-1\right)\,a^2
+\frac{\pi
\delta}{\alpha}+\frac{\pi^2\delta}{2}+\frac{l}{c_\alpha}-\frac{l^2}{2}.
\end{equation}

Combining (\ref{4.2}) and (\ref{4.4}) and simplifying, we find
that ${\rm{osc}}(u_J)=H_\alpha$, where
$H_\alpha=H_\alpha(a,l,\delta)$ is given by %

\begin{equation} \label{4.5}  %
H_\alpha=\frac{1}{2}\left(\frac{(1-lc_\alpha)^2}{1+\delta}-1\right)\,a^2+\frac{l}{1+\alpha\pi}a+\pi l+\frac{\pi^2\delta}{2}-\frac{l^2}{2}. %
\end{equation} %

Next we fix $l\in (0,\pi)$, $\alpha>0$, $\delta\ge 0$ and treat
$H_\alpha$ as a function $H_\alpha(a)$ of $a\in [0,\pi-l]$.  Since
$(1-lc_\alpha)^2/(1+\delta)<1$, it follows from (\ref{4.5}) that
$H_\alpha(a)$ is concave and takes its maximum at the
point %
\begin{equation} \label{4.6} %
a_0=\frac{(1+\delta)l}{(1+\alpha\pi)(\delta+2lc_\alpha-l^2c_\alpha^2)}.
\end{equation} %

Now, the question is whether or not the point $a_0$ given in
(\ref{4.6}) belongs to the interval $[0,\pi-l]$. Notice that the
function $\delta+2lc_\alpha-l^2c_\alpha^2$ in the denominator of
(\ref{4.6}) increases when $\alpha$ varies from $0$ to $\infty$.
This implies that, if $\delta$ and $l$ are fixed, then $a_0$,
defined by (\ref{4.6}) and considered as a function
$a_0=a_0(\alpha)$ of $\alpha$, decreases from $(1+\delta)l/\delta$
to $0$ if $\delta>0$ or from $\infty$ to $0$ if $\delta=0$, when
$\alpha$ varies from $0$ to $\infty$. Therefore if $(1+\delta)l\le
\delta(\pi-l)$, then $a_0\in [0,\pi-l]$ for all $\alpha>0$ and if
$(1+\delta)>\delta(\pi-l)$, then there is a unique $\alpha_0>0$
such that $a_0\in [0,\pi-l]$ for $\alpha\ge \alpha_0$ and
$a_0>\pi-l$ for $0<\alpha<\alpha_0$. The value of
$\alpha_0$ is the positive root of the quadratic equation %
\begin{equation} \label{4.7} %
(\pi^2\delta+2\pi l-l^2)\alpha^2+\frac{2(\pi \delta+l)(\pi-l)-\pi
l(1+\delta)}{\pi-l}\,\alpha+\delta-\frac{(1+\delta)l}{\pi-l}=0,
\end{equation} %
which is
\begin{equation} \label{4.7.1}
\alpha_0=
\frac{-2\pi^2\delta+3\pi \delta l-\pi l+2l^2+l\sqrt{(1+\delta
)(\pi^2\delta +9\pi^2-16\pi l+8l^2)}}{2(\pi^3 \delta -\pi^2 \delta
l+2\pi^2 l-3\pi l^2+l^3)}.
\end{equation}
When $\delta=0$ and $l=\pi/2$, we obtain
$\alpha_0=\frac{2}{\sqrt{3}\pi}$ that is the transitional value of
$\alpha$ found in Proposition~3.3 in \cite{LM}.

Thus, we can define the following parameters: %
\begin{equation}  \label{4.8} %
\alpha_g=\left\{\begin{array}{ll}\alpha_0 &{\mbox{if
$(1+\delta)l>\delta(\pi-l$),}}\\ %
0 & {\mbox{if $(1+\delta)l\le\delta(\pi-l)$}} %
\end{array} \right.
\end{equation} %
and
\begin{equation}  \label{4.9}
a_g=\left\{\begin{array}{ll} \pi-l &{\mbox{if
$(1+\delta)l>\delta(\pi-l$),}}\\ %
a_0 & {\mbox{if $(1+\delta)l\le \delta(\pi-l)$.}} %
\end{array} \right.
\end{equation} %

Figure~2, which illustrates possible situations, contains graphs
of ${\rm osc} (u_J)=H_\alpha(a)$ considered as a function of the
parameter $a$ for some fixed values of $l$, $\delta$ and $\alpha$.
These graphs show, in particular, that when $\alpha$ varies from
$0$ to $\infty$, the central point of the segment providing the
maximal oscillation for $u_J$ moves toward the center of the rod
from the right. For a special case when $\delta=0$, $l=\pi/2$,
this behavior had been already observed in Proposition~3.3 in
\cite{LM}.

\begin{figure}
\includegraphics[width=320pt]{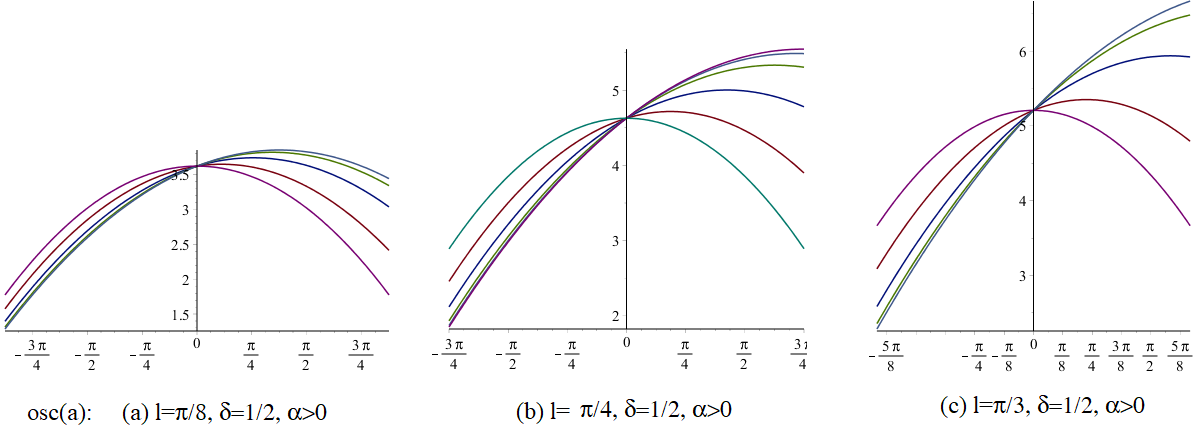}
\caption{Graphs of ${\rm osc}(u_J)=H_\alpha(a)$ as a function of
$a$ for different choices of $l$, $\delta$, and $\alpha$.}
\end{figure}

Summarizing the results of this section and evaluating the value
of the oscillation ${\rm{osc}}(u_J)=H_\alpha(a,l,\delta)$ at the
point $a=a_g$ defined in (\ref{4.9}), we obtain the following.

\begin{lemma}  %
Let $l\in (0,\pi)$, $a\in [0,\pi-l]$, $\alpha>0$, and $\delta\ge
0$ be fixed and
let $u_J$ be defined as in (\ref{4.1}). Then %
\begin{equation} \label{4.10}  %
\max_{a\in [0,\pi]}{\rm{osc}}(u_J)=\Theta_\alpha(l,\delta), %
\end{equation} %
where
\begin{equation} \label{4.11}
\Theta_\alpha=\frac{1}{2}\frac{(1+\delta)l^2}{(\pi^2\delta+2\pi l
-l^2)\alpha^2+2(\pi \delta+l)\alpha+\delta}+\pi l
+\frac{\pi^2\delta}{2}-\frac{l^2}{2},
\end{equation} %
when $(1+\delta)l<\delta(\pi-l)$ and $\alpha>0$ or $(1+\delta)l\ge
\delta(\pi-l)$ and $\alpha>\alpha_g$, and
\begin{equation}  \label{4.12}  
\Theta_\alpha=\frac{1}{2}\frac{(l^2c_\alpha^2-2lc_\alpha-\delta)(\pi-l)^2}{1+\delta}+\frac{l(\pi-l)}{1+\alpha\pi}
+\pi l+\frac{\pi^2\delta}{2}-\frac{l^2}{2}, 
\end{equation} 
when $(1+\delta)l\ge\delta(\pi-l)$ and $0<\alpha\le \alpha_g$.

Furthermore, if $l\in (0,\pi)$, $\alpha>\alpha_g$, and $\delta\ge
0$ are fixed, then ${\rm{osc}}(u_J)$ considered as a function of
$a$ increases when $a$ varies from $0$ to $a_0$ defined by
(\ref{4.6}) and ${\rm{osc}}(u_J)$ decreases when $a$ varies from
$a_0$ to $\pi-l$. If $0<\alpha\le \alpha_g$, then
${\rm{osc}}(u_J)$ increases when $a$ varies from $0$ to $\pi-l$. %
\end{lemma} %


 \section{Approximation by step functions and continuity.} %
 In this section,  we collect auxiliary
 results on convergence of sequences of solutions $u_{f_n}$ of the Robin problem,
 on approximation of $u_f$ by sequences $u_{f_n}$ with piece-wise constant functions $f_n$
 and on continuity of $u_{f_n}$ as a function of the parameters of approximants.
 These results will be used to prove our main theorems.
 We start with the following convergence lemma.

\begin{lemma}\label{T1L0}
If $f_n\in \mathcal{F}(m,M,s)$, $n=1,2,\dots$, and $f_n\to f$
($n\to\infty$) a.e. on $[-\pi,\pi]$, then $u_{f_n}\to u_f$
uniformly on $[-\pi,\pi]$ and for some subsequence,
\begin{equation}\label{5.1}
\lim_{k\to \infty} {\rm osc}(u_{f_{n_k}})={\rm osc}(u_f).
\end{equation}
\end{lemma}


\noindent \emph{Proof.} By the Dominated Convergence Theorem,
$f_n$ converges to $f$ in $L^1$. By (\ref{1.2}), for every $x\in
[-\pi,\pi]$,
\begin{eqnarray}\label{5.2}
|u_{f_n}(x)-u_f(x)|&=&\left |\int_{-\pi}^\pi G(x,y)(f_n(y)-f(y))dy\right | \\
&\leq & \max_{[-\pi,\pi]^2} G \;\;\|f_n-f\|_{1}. \nonumber
\end{eqnarray}
So,  $u_{f_n}$  converges uniformly to $u_f$ on $[-\pi,\pi]$.

For each $n\in\mathbb N$, let $x_n$, $\tilde{x}_n$ be points of
minimum and maximum of $u_{f_n}$, respectively. Choose a
subsequence such that $x_{n_k}\to x_o$ and $\tilde{x}_{n_k}\to
\tilde{x}_o$. By a standard property of uniform convergence (see
e.g. \cite[Problems 3.1.18, 3.1.23]{KN}),
$$
u_{f_{n_k}}(x_{n_k})\to
u_f(x_o),\;\;\;\;\;\;u_{f_{n_k}}(\tilde{x}_{n_k})\to
u_f(\tilde{x}_o),\;\;\;\;\;\;(k\to\infty).
$$
Therefore, for every $x\in [-\pi,\pi]$,
$$
u_f(x)=\lim_k u_{f_{n_k}}(x)\geq \lim_k
u_{f_{n_k}}(x_{n_k})=u_f(x_o).
$$
So $u_f(x_o)=\min u_f$. Similarly  $u_f(\tilde{x}_o)=\max u_f$. It
follows that
$$
{\rm
osc}(u_{f_{n_k}})=u_{f_{n_k}}(\tilde{x}_{n_k})-u_{f_{n_k}}(x_{n_k})\to
u_f(\tilde{x}_o)-u_f(x_o)={\rm osc}(u_f).
$$
The proof is complete. \hfill $\Box$

\medskip

For $n,k\in \mathbb{N}$, $1\le k\le n$, let $I_{n,k}=[-\pi+2\pi
(k-1)/n,-\pi+2\pi k/n]$. 
Thus, the intervals $I_{n,1},\ldots,I_{n,n}$ constitute a
partition of the interval $[-\pi,\pi]$ into $n$ subintervals of
equal length. We need the following approximation result.

\begin{lemma} %
Let $f\in \mathcal{F}(m,M,s)$ and $\alpha>0$. Then for every
$n,k\in \mathbb{N}$, $1\le k\le n$, there are constants
$c_{n,k}$, $m\le c_{n,k}\le M$,  such that $f_n=\sum_{k=1}^{n} c_{n,k} \chi_{I_{n,k}}$ satisfies the following: %
\begin{enumerate} %
\item[1)]  $m\le f_n\le f_{n+1}\le M$, $f_n\to f$ a.e. on
$[-\pi,\pi]$, $\|f_n\|_{L^1}\to \|f\|_{L^1}=2\pi s$.

\item[2)] $u_{f_n}(x)\to u_f(x)$ uniformly on $[-\pi,\pi]$.   %
\end{enumerate} %

\end{lemma} %

\noindent \emph{Proof.}
  Let $\varepsilon >0$. By a standard approximation result (see
  e.g. \cite[Theorem 3.14]{Ru}),   there exists a  continuous function  $g$
on $[-\pi,\pi]$ such that
  $\|f-g\|_{L^1}<\varepsilon/2$. We can demand that $m\leq g\leq M$.
  For $n\in \mathbb N$ and $k=1,2,\dots,n$, set
$$
c_{n,k}=\min\{g(x):x\in
I_{n,k}\}\;\;\;\;\;\hbox{and}\;\;\;\;\;
f_n=\sum_{k=1}^{n} c_{n,k} \chi_{I_{n,k}}.
$$
Since $f_n\leq g$ and $g$ is Riemann integrable,
$$
\|f_n-g\|_{L^1}=\int_{[-\pi,\pi]}(g-f_n)\to 0 \;\;(n\to\infty).
$$
Hence  $f_n\to g$ in $L^1$. So  $\|f_n-g\|_{L^1}<\varepsilon/2$
for all sufficiently large $n$. It follows that $f_n\to f$ in
$L^1$. Moreover, $(f_n)$ is an increasing sequence of functions.
Therefore, $f_n\to f$ a.e. on $[-\pi,\pi]$. The other assertions
come easily from Lemma \ref{T1L0}. \hfill $\Box$

Let $n\in \mathbb{N}$, $m,M,s\in \mathbb{R}$, $0\le m \le s\le M$,
and $\alpha>0$ be fixed. Let $K_n=K_n(m,M,s)$ denote the compact
set of points
$$
(\overline{t},\overline{c})=(t_1,t_2,\ldots,t_n,c_1,\ldots,c_n)\in
\mathbb{R}^{2n}
$$
such that $-\pi=t_0\le t_1\le \ldots\le t_n=\pi$,
$0\le c_k\le M-m$ and
$$
\sum_{k=1}^n c_k(t_k-t_{k-1})=2\pi (s-m).
$$
For $(\overline{t},\overline{c})\in K_n$ and $x\in [-\pi,\pi]$,
consider the function
$
f_{\overline{t},\overline{c}}\in
\mathcal{F}(m,M,s)$ defined as
$$
f_{\overline{t},\overline{c}}(x)=m+\sum_{k=1}^n c_k
\chi_{I_k}(x), \;\;\hbox{where}\;\; I_k=[t_{k-1},t_k].
$$
 Then the solution
$u_{f_{\overline{t},\overline{c}}}(x)$ of the Robin problem is a
linear combination of elementary functions as in the equation
(\ref{2.1}). This immediately implies the following continuity
lemma.

\begin{lemma} %
Let $u_{f_{\overline{t},\overline{c}}}(x)$ be the solution to the
Robin problem considered as a function of $x\in [-\pi,\pi]$ and
$(\overline{t},\overline{c})\in K_n$. Then
$u_{f_{\overline{t},\overline{c}}}(x)$ is continuous on the
compact set $[-\pi,\pi]\times K_n$.

In particular, if $E_1$ is a compact subset of $[-\pi,\pi]$, $E_2$
is a compact subset of $K_n$, then there are points in $E_1\times
E_2$, where $u_{f_{\overline{t},\overline{c}}}(x)$ achieves its
minimum and maximum on $E_1\times E_2$.
\end{lemma} %

 \section{Main proofs.}
Now we are ready to present proofs of our main results. We start
with the proof of Theorem~2.

\medskip

\noindent %
{\bf Proof of Theorem 2.} Let $f\in \mathcal{F}(m,M,s)$. First,
we prove the right inequality in (\ref{1.6}). It follows from
Lemma~5 that, for every $\varepsilon>0$,
there exists a piece-wise constant function
$f_{n,\overline{c}'}=m+\sum_{k=1}^n c_{n,k}' \chi_{I_{n,k}}\in
\mathcal{F}(m,M,s_1)$, where
$\overline{c}'=(c'_{n,1},\ldots,c'_{n,n})$, such that $0\le
c_{n,k}'\le M-m$, $|s_1-s|<\varepsilon$ and,  for the given $x_0\in
[-\pi,\pi]$,%
\begin{equation}  \label{6.1} %
|u_{f_{n,\overline{c}'}}(x_0)-u_f(x_0)|<\varepsilon.
\end{equation}  %

We note here that the integer $n$ in the definition of the
function $f_{n,\overline{c}'}$  can be chosen as large as we need
for our proof. Indeed, for any integer $j\ge 1$, consider the
partition of $[-\pi,\pi]$ into the intervals $I_{nj,s}$, $1\le
s\le nj$. Then consider
$\overline{c}''=(c_{nj,1}'',\ldots,c_{nj,nj}'')$ such that
$c_{nj,s}''=c_{n,k}'$ if $I_{nj,s}\subset I_{n,k}$. With these
notations, we have $f_{nj,\overline{c}''}=f_{n,\overline{c}'}$
a.e. on $[-\pi.\pi]$.

 Next,
for $n$ sufficiently large, we consider
$f_{\overline{c}}(x_0)=m+\sum_{k=1}^n c_k \chi_{I_{n,k}}(x_0)$ 
and $u_{f_{\overline{c}}}(x_0)$ as functions of
$\overline{c}=(c_1,\ldots,c_n)$, assuming that $\overline{c}$
varies over the compact set defined by the following conditions:
\begin{equation} \label{6.2} 
0\le c_k\le M-m, \quad \frac{1}{n}\sum_{k=1}^n c_k=s_1-m. %
\end{equation} %

It follows from the continuity Lemma~6, that there is
$\overline{c}^*=(c_1^*,\ldots,c_n^*)$, satisfying  conditions
(\ref{6.2}), such that
\begin{equation} \label{6.3} %
\max u_{f_{\overline{c}}}(x_0)\le u_{f_{\overline{c}^*}}(x_0),
\end{equation}  %
where the maximum is taken over all $\overline{c}$ satisfying
conditions (\ref{6.2}).

Let $a^*=a_m(\frac{\pi}{n})$ denote the center of the interval $I$
of length $\frac{2\pi}{n}$, which provides the maximum to the
function $u_{\chi_I}(x)$ at the point $x=x_0$. We recall that
$a^*=a_m(\frac{\pi}{n})$, given by equation (\ref{2.3}), is such
that $0\le a^*\le x_0\le a^*+\frac{\pi}{n}$.


Suppose that the centers of the intervals $I_{n,k}$ and
$I_{n,k+1}$ lie in the interval $[-\pi,a^*]$. Under these
conditions, we claim that either $c_k^*=0$ or $c_{k+1}^*=M-m$.
Indeed, if $c_k^*>0$ and $c_{k+1}^*<M-m$ then, for all
sufficiently small $\delta>0$, $\widetilde{f}=f_{\overline{c}^*}-\delta
\chi_{I_{n,k}}+\delta\chi_{I_{n,k+1}}\in \mathcal{F}(m,M,s_1)$.
 Furthermore, it follows from Lemma~1 that
 $u_{\chi_{I_{n,k}}}(x_0)<u_{\chi_{I_{n,k+1}}}(x_0)$. The latter together with the
 additivity property (\ref{1.4}) imply that
 $u_{\widetilde{f}}(x_0)>u_{f_{\overline{c}^*}}(x_0)$,
 contradicting the inequality (\ref{6.3}).
 A similar argument shows that if the centers of the intervals $I_{n,k}$ and
$I_{n,k+1}$ lie in the interval $[a^*,\pi]$ then either
$c_k^*=M-m$ or $c_{k+1}^*=0$.

Using these observations, we conclude that, if the integer $n$ is
large enough, then there are integers $k_1,k_2$, $1\le k_1<k_2\le
n$, such that $c_j^*=0$ for $1\le j<k_1$ and $k_2<j\le n$,
$c_j^*=M-m$ for $k_1< j< k_2$, and $0\le c_j^*<M-m$, for $j=k_1$
and $j=k_2$.

Since $f_{n,\overline{c}'}\in \mathcal{F}(m,M,s_1)$, it follows
from (\ref{6.3}) and the additivity property (\ref{1.4}) that %
\begin{equation}  \label{6.4} %
u_{f_{n,\overline{c}'}}(x_0)\le
u_{f_{\overline{c}^*}}(x_0)<u_{\widehat{f}}(x_0),
\end{equation} %
where $\widehat{f}=m+(M-m)\chi_{\widehat{I}}$ and
$\widehat{I}=[-\pi+2\pi k_1/n,-\pi+2\pi k_2/n]$.

Let $l_1=\pi(s_1-m)/(M-m)$ and let $\widehat{l}$ denote the half
length of the interval $\widehat{I}$. Since $f^*\in F(m,M,s_1)$,
it follows from our construction of $\widehat{I}$ that
\begin{equation} \label{6.5} %
l_1\le \widehat{l}\le l_1+\frac{2\pi}{n}.
\end{equation} %

It follows from Lemma~1 and equation (\ref{6.5}) that %
\begin{equation} \label{6.6} %
u_{\chi_{\widehat{I}}}(x_0)\le \nu_\alpha(x_0,\widehat{l})\le
\nu_\alpha(x_0,l_1)+O(1/n),
\end{equation} %
where $\nu_\alpha(x,l)$ is defined by (\ref{2.4}). The latter
inequality together with equation (\ref{2.2}) and the additivity
property (\ref{1.4}) implies the following: %
\begin{equation} \label{6.7} %
u_{\widehat{f}}(x_0)\le \eta_\alpha(x_0)
m+(M-m)\nu_\alpha(x_0,l_1)+O(1/n).
\end{equation} %

Combining (\ref{6.1}), (\ref{6.4}) and (\ref{6.7}), we conclude
that %
\begin{equation} \label{6.8} %
u_{f}(x_0)\le u_{f_{n,\overline{c}'}}(x_0)+\varepsilon \le
\eta_\alpha(x_0) m+(M-m)\nu_\alpha(x_0,l)+\delta(\varepsilon,n),
\end{equation} %
where $\delta(\varepsilon,n)\to 0$ when $\varepsilon\to 0$ and
$n\to \infty$. Since $\varepsilon>0$ in equation (\ref{6.1})  can
be chosen arbitrarily small and the integer $n$ can be chosen
arbitrarily large, (\ref{6.8}) implies the second inequality in
(\ref{1.6}).

To prove uniqueness, we assume that $f$ does not coincide with
$f_0^+=m+(M-m)\chi_{I(a_m,l)}$ on a set of positive measure.
Since $f\in L^1$, there are ``density points'' $x_1$ and $x_2$,
such that 
$f_0^+(x_1)=M$ and
$f(x) <M-\varepsilon$ on a subset $E_1$ of some small interval
$I_1$ centered at $x_1$, $f_0^+(x_2)=m$ and $f(x)
>m+\varepsilon$ on a subset $E_2$ of some small interval $I_2$ centered
at $x_2$. Let $E=E_1\cap (E_2+(x_1-x_2))$. Since $x_1$ and $x_2$
are density points for $f\in L^1$, the one-dimensional Lebesgue measure of $E$ is strictly positive. We may assume without loss of generality that
$E$ and $E-(x_1-x_2)$ either both lie on the interval $[-\pi,a_m]$
or both lie on the interval $[a_m,\pi]$. Under these assumptions,
it follows from the monotonicity properties of Lemma~1 that
$u_{\chi_E}(x_0)>u_{\chi_{E+(x_1-x_2)}}(x_0)$. Therefore,
replacing $f$ by $\widetilde{f}=f-\varepsilon
\chi_{E-(x_1-x_2)}+\varepsilon \chi_E$ with $\varepsilon>0$ small
enough, we obtain a function $\widetilde{f}\in \mathcal{F}(m,M,s)$
such that $u_f(x_0)<u_{\widetilde{f}}(x_0)$. Since
$\widetilde{f}\in \mathcal{F}(m,M,s)$, the right inequality in
(\ref{1.6}) holds true for $u_{\widetilde{f}}$. Therefore, for the
function $u_f$ the right inequality in (\ref{1.6}) holds with the
sign of strict inequality.

\medskip

To prove the left inequality in (\ref{1.6}), we assume that $f\in
\mathcal{F}(m,M,s)$ and consider the function $f^-=M+m-f$. We have
$m\le f^-\le M$ and $\|f^-\|_{L^1}=2\pi s^-$ with $s^-=M+m-s$,
$m<s^-<M$. Thus, $f^-\in \mathcal{F}(m,M,s^-)$ and therefore, by
our proof above,
\begin{equation} \label{6.9} %
u_{f^-}(x_0)\le \eta_\alpha(x_0)m+(M-m)\nu_\alpha(x_0,l^-),
\end{equation} %
where $l^-=\pi(s^--m)/(M-m)=\pi-l$.
Since $f+f^-=M+m$, we have %
\begin{equation} \label{6.10} %
u_f(x)+u_{f^-}(x)=\int_{-\pi}^\pi
G(x,y)(f(y)+f^-(y))\,dy=(M+m)\eta_\alpha(x).
\end{equation} %
Combining equations (\ref{6.9}) and (\ref{6.10}), we obtain the
left inequality in (\ref{1.6}).

Furthermore,  (\ref{1.6}) holds with the equality sign in the left
inequality if and only if (\ref{6.9}) holds with the sign of
equality. We proved above that the latter holds if and only if
$f^-=m+(M-m)\chi_{I(a_m^-,l^-)}$ a.e. on $[-\pi,\pi]$, which
implies that the sign of equality in the left inequality in
(\ref{1.7}) occurs if and only if $f=f_0^-$ a.e. on $[-\pi,\pi]$. %
 \hfill $\Box$

 \bigskip

 The structure of proofs of Theorems 1 and 3 presented below is the same as in the proof of
 Theorem~2. Basically, to prove these theorems, we use the same arguments
 as in the proof of Theorem~2 with minor
 changes. Below, we sketch
 these proofs emphasizing these minor changes.

 \medskip

\noindent %
 {\bf Proof of Theorem 3.}
 First, given $\varepsilon>0$, we approximate $f\in \mathcal{F}(m,M,s)$
 with $f_{n,\overline{c}'}\in \mathcal{F}(m,M,s_1)$ such that
 $|s_1-s|<\varepsilon$ and
\begin{equation}  \label{6.11} %
|(u_{f_{n,\overline{c}'}}(x_0)-u_{f_{n,\overline{c}'}}(-\pi))
-(u_f(x_0)-u_f(-\pi))|<\varepsilon.
\end{equation}
Then, using the continuity Lemma~6, we find $f^*\in \mathcal{F}(m,M,s_1)$ such that  %
\begin{equation} \label{6.12} %
\max \{u_{f_{\overline{c}}}(x_0)-u_{f_{\overline{c}}}(-\pi)\}\le
u_{f_{\overline{c}^*}}(x_0)-u_{f_{\overline{c}^*}}(-\pi),
\end{equation}  %
where the maximum is taken over all $\overline{c}$ satisfying
conditions (\ref{6.2}).

Next, performing tricks with the intervals $I_{n,k}$ and
$I_{n,k+1}$ as we did in the proof of Theorem~2 and using Lemma~2,
we obtain a function $\widehat{f}=m+(M-m)\chi_{\widehat{I}}$,
where $\widehat{I}=[-\pi+2\pi k_1/n,-\pi+2\pi k_2/n]$ with
appropriate $k_1$ and $k_2$, such that
the following holds: %
\begin{eqnarray}\label{6.13}
u_{f_{\overline{c}^*}}(x_0)&-&u_{f_{\overline{c}^*}}(-\pi)<u_{\widehat{f}}(x_0)-u_{\widehat{f}}(-\pi)
\le \\ &{ }&
m(\eta_\alpha(x_0)-\eta_\alpha(-\pi))+(M-m)\tau_\alpha(x_0,l_1)+O(1/n).
\nonumber
\end{eqnarray} %
Since $\varepsilon>0$ in (\ref{6.11}) can be taken arbitrarily
small and the integer $n$ in (\ref{6.13}) can be taken arbitrarily
large, combining (\ref{6.11}), (\ref{6.12}) and (\ref{6.13}), we
obtain the inequality in (\ref{1.7}).

The proof of the uniqueness statement of Theorem~3, is almost
identical with the uniqueness proof of Theorem~2. Namely, given
$f\in \mathcal{F}(m,M,s)$, we use our ``two density points
argument'', to construct a function $\widetilde{f}=f-\varepsilon
\chi_{E-(x_1-x_2)}+\varepsilon \chi_E \in \mathcal{F}(m,M,s)$,
with $\varepsilon>0$ small enough,  such that
$u_f(x_0)-u_f(-\pi)<u_{\widetilde{f}}(x_0)-u_{\widetilde{f}}(-\pi)$.
Since $\widetilde{f}\in \mathcal{F}(m,M,s)$, the inequality in
(\ref{1.7}) holds true for $u_{\widetilde{f}}$. Therefore, for the
function $u_f$,  (\ref{1.7}) holds with the sign of strict
inequality. %
 \hfill $\Box$

 \bigskip

\noindent %
{\bf Proof of Theorem 1.} %
 It follows from Lemma~5 that, for given $f\in F(m,M,s)$ and
 $\varepsilon>0$ arbitrarily small, there exists  a piece-wise constant function
 $f_{n,\overline{c}'}=m+\sum_{k=1}^n c_{n,k}' \chi_{I_{n,k}}\in
\mathcal{F}(m,M,s_1)$, where
$\overline{c}'=(c'_{n,1},\ldots,c'_{n,n})$, such that $0\le
c_{n,k}'\le M-m$, $|s_1-s|<\varepsilon$ and  %
\begin{equation}  \label{6.14} %
|{\rm osc}(u_{f_{n,\overline{c}'}})-{\rm osc}(u_f)|<\varepsilon.
\end{equation}  %
Here, we assume once more, that the integer $n$ is chosen as large
as we need for our proof.

Then, using the continuity Lemma~6, we find $f_{\overline{c}^*}=m+\sum_{k=1}^n c_{n,k}^* \chi_{I_{n,k}}\in \mathcal{F}(m,M,s_1)$ such that  %
\begin{equation} \label{6.15} %
\max \{{\rm osc}(u_{f_{\overline{c}'}})\}\le {\rm
osc}(u_{f_{\overline{c}^*}}),
\end{equation}  %
where the maximum is taken over all $\overline{c}$ satisfying
conditions (\ref{6.2}). Since $u_{f_{\overline{c}^*}}$ is concave
on $[-\pi,\pi]$, we may assume without loss of generality that
${\rm
osc}(u_{f_{\overline{c}^*}})=u_{f_{\overline{c}^*}}(x_0)-u_{f_{\overline{c}^*}}(-\pi)$
for some $x_0\in (-\pi,\pi]$. If this is not the case, we replace
$f_{\overline{c}^*}(x)$ with $f_{\overline{c}^*}(-x)$.

Next, we consider the temperature gap
$E(a)=u_{\chi_{I(a,l_n)}}(x_0)-u_{\chi_{I(a,l_n)}}(-\pi)$ for the
interval $I(a,l_n)$ centered at $a$ with half length
$l_n=\pi/n$. As we have shown in Section~3, there is a unique
point $a^*=a_e$ with $a_e$ given by equation (\ref{3.4}), where
$E(a)$ achieves its maximum $\tau_\alpha$ given by equation
(\ref{3.5}).

If the centers of the intervals $I_{n,k}$ and $I_{n,k+1}$ both lie
in the interval $[-\pi,a^*]$ or both lie in the interval
$[a^*,\pi]$ then arguing as in the proof of Theorem~2 and using
the monotonicity properties of Lemma~2, we conclude  that if
$f_{\overline{c}^*}$ is maximal in the sense of equation
(\ref{6.15}),
then either $c_k^*=0$ or $c_{k+1}^*=M-m$.

The latter implies that there exists an interval
$\widehat{I}=[-\pi+2\pi k_1/n,-\pi+2\pi k_2/n]$ with half
length $\widehat{l}=\pi(k_2-k_1)/n$ such that
$l_1<\widehat{l}<l_1+\pi/n$, where $l_1=\pi(s_1-m)/(M-m)$, and
such that for $\widehat{f}=m+(M-m)\chi_{\widehat{I}}$ we have
\begin{equation} \label{6.16} %
u_{f_{\overline{c}^*}}(x_0)-u_{f_{\overline{c}^*}}(-\pi)\le
u_{\widehat{f}}(x_0)-u_{\widehat{f}}(-\pi)\le {\rm
osc}(u_{\widehat{f}}).
\end{equation}  %

Now, it follows from equations (\ref{4.10}), (\ref{4.11}),
(\ref{4.12}) of Lemma~3 that %
\begin{equation} \label{6.17} 
{\rm osc}(u_{\widehat{f}})\le
(M-m)\Theta_\alpha(\widehat{l},\delta)=(M-m)\Theta_\alpha(l,\delta)+\beta(\varepsilon,n),
\end{equation}
where $\delta=m/(M-m)$ and $\beta(\varepsilon,n)\to 0$ when
$\varepsilon\to 0$ and $n\to \infty$.

Finally, combining equations (\ref{6.14}) --(\ref{6.17}), we
obtain the inequality in (\ref{1.5}).

To prove the uniqueness statement of Theorem~1,
we use once more the ``two density points argument'' as in the
proofs of Theorems~2 and 3, to construct a function
$\widetilde{f}=f-\varepsilon \chi_{E-(x_1-x_2)}+\varepsilon \chi_E
\in \mathcal{F}(m,M,s)$, with $\varepsilon>0$ small enough,  such
that ${\rm osc}(u_f)<{\rm osc}(u_{\widetilde{f}})$. Since
$\widetilde{f}\in \mathcal{F}(m,M,s)$, the inequality in
(\ref{1.5}) holds true for $u_{\widetilde{f}}$. Therefore, for the
function $u_f$, (\ref{1.5}) holds with the sign of strict
inequality. %

\hfill $\Box$

\bigskip

\section{Temperature gap in higher dimensions} %
One can consider a variety of higher dimensional analogs of
Problem 1. Here we present two such problems on the temperature
gap in pipes $P_L$, which are cylindrical domains
$P_L=\{(x_1,x_2,x_3)\in \mathbb{R}^3:\,x_1^2+x_2^2<1,\,|x_3|<L\}$,
$L>0$. Let $S_L$ denote the cylindrical boundary of $P_L$ and let
$D_L\ni (0,0,L)$ and  $D_{-L}\ni (0,0,-L)$ denote the boundary
disks of $P_L$. By $\frac{\partial}{\partial {\rm n}}$ we denote
the outward normal derivative on $\partial P_L$ (for the boundary
points where it is defined).

\begin{problem}
        Let $E$ be a compact subset of $P_L$ of given volume $V$,
        $0<V<2\pi L$,  and let $\alpha>0$. Suppose that
        $u_E$ is a bounded solution to the Poisson equation
        $$ 
        -\Delta u=\chi_E
        \;\;\;\;\hbox{in}\;\;\; P_L, %
        $$ 
        with mixed Neumann-Robin boundary conditions %
        $$ 
        \frac{\partial u}{\partial {\rm n}}=0 \;\;\;\;\hbox{on}\;\;\; S_L, %
        $$ 
        $$ 
        \frac{\partial u}{\partial {\rm n}}+\alpha u=0
        \;\;\;\;\hbox{on}\;\;\;D_{\pm L}.
        $$ 
Find $\max {\rm osc}(u_E)$ over all open sets $E\subset P_L$ of
volume $V$ and identify sets $E^*$ providing this maximum.

\end{problem}

\begin{problem}
Let $\Omega$ be a compact subset of $S_L$ of given area $A$,
$0<A<4\pi L$  and let $\alpha>0$. Suppose that $v_\Omega$ is a
bounded solution to the Laplace equation
$$ 
\Delta v=0 \;\;\;\;\hbox{in}\;\;\; P_L,
$$ 
with mixed Dirichlet-Robin boundary conditions
$$ 
v=\chi_\Omega \;\;\;\;\hbox{on}\;\;\; S_L,
$$ 
$$ 
\frac{\partial v}{\partial {\rm n}}+\alpha v=0
\;\;\;\;\hbox{on}\;\;\;D_{\pm L}.
$$ 

Find $\max {\rm osc}(v_\Omega)$ over all open sets $\Omega\subset
S_L$ of area $A$ and identify sets $\Omega^*$ providing this
maximum.
\end{problem}

We assume that solutions of Problems 2 and 3 exist and are
regular; for the existence and regularity of elliptic problems
with Robin boundary conditions, we refer to \cite{Ni2011}  and
references therein.
 Problems~2 and 3 look challenging. Similar
problems to find $\max u_E$ or  $\max v_\Omega$ on $P_L$ instead
of ${\rm osc} (u_E)$ in Problem~2 or ${\rm osc}(v_\Omega)$ in
Problem~3 can be more tractable. We expect that the symmetric
configurations will provide the corresponding maxima and therefore
the symmetrization methods due to Talenti and Baernstein can be
applied to solve these problems.

One can also consider analogs of Problems~2 and 3 in cylindrical
domains in $\mathbb{R}^n$ of any dimension $n\ge 2$.

\bigskip %

{\bf Acknowledgments.} We thank Dr. G. Sakellaris and Dr. K.
Yamazaki for helpful discussions.


\begin{thebibliography}{aaaa}  %


 \bibitem{A} {\sc S. Abramovich}, \emph{Monotonicity of eigenvalues under
 symmetrization}.
    SIAM J. Appl. Math., \textbf{28} (1975), no. 2, 350--361.

    \bibitem{B}
    {\sc A. Baernstein II,} \emph{Symmetrization in analysis.}
    With David Drasin and Richard S. Laugesen. With a foreword by
    Walter Hayman. New Mathematical Monographs, 36.
    Cambridge University Press, Cambridge, 2019. 


 \bibitem{BS} {\sc D. Betsakos and A. Yu. Solynin,}  \emph{Heating long pipes.} Anal. Math. Phys. \textbf{11} (2021), no. 1,
 Paper No. 40, 35 pp. 

  \bibitem{KN} {\sc W.J. Kaczor and M.T. Nowak}, \emph{Problems in Mathematical Analysis. II. Continuity and differentiation}.
  Translated from the 1998 Polish original, revised and augmented by the authors.
  Student Mathematical Library, 12. American Mathematical Society, Providence, RI, 2001


    \bibitem{LM} {\sc J. J. Langford and P. McDonald}, \emph{Extremizing temperature functions of rods with Robin boundary conditions}.
    Ann. Fenn. Math. \textbf{47} (2022), no. 2, 759-775.

\bibitem{Ni2011} {\sc R. Nittka}, Regularity of solutions of linear second order elliptic
and parabolic boundary value problems on Lipschitz domains. J.
Differential Equations 251 (2011), no. 4--5, 860--880.

\bibitem{PS} {\sc G. P\'{o}lya and G.~Szeg\"{o}}, \emph{Isoperimetric Inequalities in Mathematical Physics}.
    Princeton University Press, Princeton, N.J., 1951.

\bibitem{Ru} {\sc W. Rudin}, \emph{Real and Complex Analysis}. Third edition. McGraw-Hill Book Co., New York, 1987.


 \bibitem{T} {\sc G. Talenti}, \emph{The art of rearranging}.  Milan J. Math., 84:1, 2016, 105--157.


\end{thebibliography}
\end{document}